\documentclass[12pt,a4paper,draft]{amsart}
\usepackage{color}
\usepackage[colorlinks=true]{hyperref}
\hypersetup{urlcolor=blue, citecolor=blue, draft=false}

\theoremstyle{remark}

\theoremstyle{remark}

\begin{document}

\begin{center}
{\bf MONOTONICITY AND BOUNDED VARIATION FOR HARMONIC MAPPINGS WITH APPLICATIONS TO SPECIAL FUNCTIONS}\\

\bigskip
JUGAL KISHORE PRAJAPAT \\
Department of Mathematics, Central University of Rajasthan, \\
Bandarsindri-305817, Kishangarh, Rajasthan, India\\

Email: jkprajapat@curaj.ac.in
\end{center}

\vspace{1cm}
\noindent
{\bf Abstract:} We investigate sufficient conditions for the close-to-convexity and univalence of planar harmonic mappings. First, we establish a unified monotonicity criterion for harmonic mappings that incorporates the ordinary and odd harmonic cases. We then develop a criterion based on the total variation of successive coefficient differences. Unlike conventional monotonicity conditions, this approach allows the coefficient sequences to exhibit non-monotone behavior. Several previously established criteria for harmonic close-to-convexity are recovered as particular consequences. Applications are developed for harmonic mappings associated with Gaussian hypergeometric functions and Legendre polynomials, with explicit restrictions on their parameters.



\section{Introduction}
Coefficient criteria for univalence form a standard topic in geometric function theory. Among the different approaches developed for this purpose, the results expressed directly in terms of its series coefficients are useful, since they often permit geometric conclusions without requiring an explicit description of the image domain. 

Let $\mathbb D=\{z\in\mathbb C:|z|<1\}$ denote the open unit disk, and let $\mathcal A$ be the class of analytic functions normalized by $\,h(z)=z+\sum_{n=2}^{\infty}a_nz^n,
\; z\in\mathbb D.$ We denote by $\mathcal{S}$, the subclass of univalent functions in $\mathcal{A}$. The study of coefficient criteria ensuring membership of $f \in \mathcal{A}$ in $\mathcal S$ has a long history. Alexander \cite{Alexander} investigated the connection between geometric properties of analytic functions and their derivatives. In subsequent works, Fej\'er, Szeg\H{o} and Robertson \cite{Fejer, Szego, Robertson} extended these results by examining sequences of real numbers with monotone properties.

A complex-valued harmonic function $f$ defined in a simply connected domain admits a canonical decomposition $f=h+\overline g,$ where $h$ and $g$ are analytic. Unlike analytic mappings, harmonic mappings need not preserve orientation, and their local univalence depends on derivatives of both $h$ and $g$. Therefore, coefficient conditions that ensure univalence in the analytic case do not automatically extend to harmonic mappings. 

For a normalized harmonic mapping
$$f(z)=z+\sum_{n=2}^{\infty}a_nz^n+
\overline{\sum_{n=2}^{\infty}b_nz^n},$$
the associated analytic family is $F_\varepsilon(z)=h(z)+\varepsilon g(z),
\; |\varepsilon|=1.$
If members of this family are close-to-convex, the corresponding harmonic mappings are close-to-convex and univalent. Recently, Sahoo and Wankhede \cite{Sahoo} established several sufficient conditions involving monotonicity of coefficients and their successive differences. They also obtained applications to Gaussian hypergeometric functions. However, many power series may fail to be monotone even when the corresponding functions possess geometric properties. This observation motivates the search for criteria that control the total variation rather than requiring successive differences to have the same sign. A related motivation comes from lacunary power series. Consider the harmonic mapping
$$f_p(z)=z+\sum_{n=1}^{\infty}a_{pn+1}z^{pn+1}
+\overline{\sum_{n=1}^{\infty}b_{pn+1}z^{pn+1}},$$
where $p=1$ or $p=2$. The choice $p=1$ corresponds to an unrestricted normalized harmonic mapping, whereas $p=2$ gives an odd harmonic mapping. Although these cases have been studied separately, their coefficients suggest that a common argument may be developed by introducing the weighted sequences $\,c_n=(pn+1)a_{pn+1}\,$ and $\,d_n=(pn+1)b_{pn+1}.$ This formulation allows the two cases to be studied through a single coefficient inequality.

For complex numbers $a, b$ and $c$ with $c \notin \{0, -1, -2, \cdots \},$ the Gaussian hypergeometric function $_2F_1(a, b; c; z)$ is defined by the series
\begin{align} \label{e1.1}
 _{2}F_{1}(a,b;c;z) = \displaystyle{\sum_{n=0}^{\infty}} \dfrac{(a)_{n}(b)_n}{(c)_{n}} \dfrac{z^{n}}{n!},
\end{align}
where $(\alpha)_n$ is the Pochhammer symbol defined as $(\alpha)_0=1$ and $(\alpha)_n=\alpha(\alpha+1) \cdots (\alpha+n-1)$ with $n \geq 1$. Also note that
\begin{equation}\label{eq:2}
{}_2F_1'(a,b;c;z)=
\frac{ab}{c}{}_2F_1(a+1,b+1;c+1;z).
\end{equation}
The Gaussian hypergeometric function contains many classical special functions as particular cases. In particular, the hypergeometric representation of the Legendre function is given by 
\begin{equation}\label{eq-4}
P_{\nu}(x)={}_2F_1\left(
-\nu,\nu+1;1;\frac{1-x}{2}
\right).
\end{equation}
Putting $\nu=-a,\; x=1-2z,$ we obtain $\,{}_2F_1(a,1-a;1;z)=P_{-a}(1-2z).$ The basic information on hypergeometric functions and their special cases may be found in \cite{Andrews}.

The coefficients of hypergeometric functions possess explicit recurrence relations, while their asymptotic behavior can be expressed through gamma-function identities. The geometric theory of hypergeometric functions has been investigated recently by several researchers (see, for example, \cite{Carlson, Hasto, Kumari, Miller, Ponnu, Ponnu2, Sugawa, Yildiz}). More recently, investigations of shifted hypergeometric functions and their geometric properties have further emphasized the importance of understanding the influence of parameters on univalence and convexity (see \cite{Kustner, Sugawa1, Sugawa2}). 

In this paper, we first use the monotonicity of weighted coefficients to provide a unified close-to-convexity result for ordinary and odd harmonic mappings. In another result, we show that the total variation of weighted coefficient differences yields a criterion that does not require monotonicity. Recall that the total variation of a sequence $\{x_n\}_{n\ge1}$ is defined by $\,V(x)=\sum_{n=1}^{\infty}|x_{n+1}-x_n|.$ The sequence is said to have bounded variation if $V(x)<\infty$. Every bounded-variation sequence converges, but such a sequence need not be monotone. In particular, the condition $\,V(x)<\infty$ permits successive differences to change sign infinitely often, provided their absolute values are summable. This observation is essential for extending monotonicity-based criteria to a wider family of harmonic mappings.

Let $\mathcal H$ denote the normalized class of harmonic mappings $f(z)=h+\overline{g}$ in $\mathbb{D}$ satisfying $\,h(0)=g(0)=0,\; h'(0)=1,\; g'(0)=0.$ Thus
$$h(z)=z+\sum_{n=2}^{\infty}a_nz^n, \quad \text{and} \quad 
\qquad g(z)=\sum_{n=2}^{\infty}b_nz^n, \qquad z \in \mathbb{D}.$$
The Jacobian of $f$ is $\,J_f(z)=|h'(z)|^2-|g'(z)|^2.$ The mapping $f$ is locally sense-preserving when $J_f(z)>0$ in $\mathbb D$. Whenever $h' \neq 0$ in $\mathbb{D}$, the analytic dilatation is defined by $\,\omega_f(z)=g'(z)/h'(z).$ Hence the sense-preserving condition becomes $\,|\omega_f(z)|<1.$

An analytic univalent function $h$ is said to be convex in $\mathbb D$ if $h(\mathbb D)$ is a convex domain. The standard analytic characterization for convex functions in $\mathbb{D}$ is $f'(0) \neq 0\,$ and $\,\text{Re}\left(1+zf''(z)/f'(z)\right)>0, \; z\in\mathbb D.$ A function $f\in \mathcal{A}$ is said to be close-to-convex, if there exist a convex function $g\in\mathcal A$ and a real number $\theta$ such that $\text{Re}\left(e^{i\theta} f'(z)/g'(z)\right)>0$. Every close-to-convex function is univalent, although the converse does not hold. The basic information on analytic convexity and close-to-convexity may be found in \cite{Duren,Goodman}.

The convex functions $\,\phi_1(z)=-\log(1-z)$ and $\,\phi_2(z)=(1/2)\log \,[(1+z)/(1-z)]$ play an important role in our arguments. We unify these functions by 
$$\phi_p(z)=\int_0^z \frac{1}{1-t^p}\,dt \qquad (z \in \mathbb{D}, \,\,p=1, 2).$$
Clearly, $\phi_p$ is convex for $p=\{1,2\}.$ 

\medskip
We shall need the following result to prove the main results:

\medskip
\noindent
{\bf Lemma 1.1\,}[Clunie–Sheil-Small \cite{Clunie}]. Let $h$ and $g$ be analytic in $\mathbb D$, with $\,|h'(0)|>|g'(0)|.$ If the analytic function $\,F_\varepsilon=h+\varepsilon g\,$ is close-to-convex in $\mathbb D$ for every $|\varepsilon|=1$, then $\,f=h+\overline g$ is a close-to-convex and univalent in $\mathbb D$.

\section{Main Results}

The following theorem combines coefficient criteria for ordinary harmonic functions and odd harmonic functions.

\medskip
\noindent
{\bf Theorem 2.1}
Let $p \in\{1, 2\}$ and let $\{a_{pn+1}\}_{n\geq0}$ and $\{b_{pn+1}\}_{n\geq 0}$ be sequences of real numbers with $a_1=1$ and $b_1=0$. Suppose that $\,\{(pn+1)a_{pn+1}\}_{n\geq0}\,$ and $\,\{(pn+1)b_{pn+1}\}_{n\geq0}\,$ are non-decreasing sequences and $\,\lim_{n\to\infty}(pn+1)
(a_{pn+1}+b_{pn+1})\leq2.$ Then the harmonic function
$$f_p(z)=z+\sum_{n=1}^{\infty}a_{pn+1}z^{pn+1}
+\overline{\sum_{n=1}^{\infty}b_{pn+1}z^{pn+1}}=h_p(z)+\overline{g_p(z)}$$
is close-to-convex and univalent in $\mathbb D$.

\medskip
\noindent
{\it Proof.} Let $\,c_n=(pn+1)a_{pn+1},\;
d_n=(pn+1)b_{pn+1},\; n\geq0,$ where $c_0=1$ and $d_0=0$.
By hypothesis, $c_n\geq1,\; d_n\geq0,$ and $\,1\leq c_n+d_n\leq2.$
Thus the sequences $a_{pn+1}$ and $b_{pn+1}$ are bounded, hence the series for $h$ and $g$ converge locally uniformly in $\mathbb D$, and $f_p=h_p+\overline{g_p}$ is harmonic in $\mathbb D$.

For every $\varepsilon\in\mathbb C$ with $|\varepsilon|=1$, consider $F_\varepsilon(z)=h_p(z)+\varepsilon g_p(z).$ Hence
$$(1-z^p)F_\varepsilon'(z)=1+\sum_{n=1}^{\infty}
\left[(c_n-c_{n-1})\right. \left.
+\varepsilon(d_n-d_{n-1})\right]z^{pn}.$$

For $N\geq1$, define
$$S_N(z)=1+\sum_{n=1}^{N}
\big[(c_n-c_{n-1})
+\varepsilon(d_n-d_{n-1})\big]z^{pn}.$$
Since the coefficient sequences $\{c_n\}_{n \ge 1}$ and $\{d_n\}_{n \ge 1}$ are non-decreasing, for $z \in \mathbb{D}$,
$$\begin{aligned}
\text{Re}S_N(z)
&\geq1-\sum_{n=1}^{N}|c_n-c_{n-1}||z|^{pn}-\sum_{n=1}^{N}|d_n-d_{n-1}||z|^{pn}\\
&\ge 1-\sum_{n=1}^{N}(c_n-c_{n-1})-\sum_{n=1}^{N}(d_n-d_{n-1})\\
&=2-(c_N+d_N)\geq0.
\end{aligned}$$
Choosing $N \to \infty$, we obtain $\,\text{Re}\{(1-z^p)F_\varepsilon'(z)\}>0,
\; z\in\mathbb D.$ Now define $\phi_p(z)=\int_0^z(1-t^p)^{-1} \,dt.$ Then $\phi_p$ is convex for $p=1,2$. Moreover
$$\text{Re}\left\{
\frac{F_\varepsilon'(z)}{\phi_p'(z)}
\right\}=\text{Re}\{(1-z^p)F_\varepsilon'(z)\}>0.$$
Hence, $F_\varepsilon$ is close-to-convex for every $|\varepsilon|=1$. Now the use of Lemma 1.1 shows that the harmonic mapping $f_p$ is close-to-convex and univalent in $\mathbb{D}$.

To prove the sharpness, we assume for every $n\ge1, \;a_{pn+1}=1/(pn+1),
\; b_{pn+1}=(1+\delta)/(pn+1).$ Thus both weighted sequences $\{c_n\}$ and $\{d_n\}$ are non-decreasing, and $\,\lim_{n\to\infty}(pn+1)
(a_{pn+1}+b_{pn+1})=2+\delta.$ The corresponding harmonic mapping is $\,f_{p,\delta}(z)=h_p(z)+\overline{g_{p,\delta}(z)},$ where
$$h_p(z)=\sum_{n=0}^{\infty}
\frac{z^{pn+1}}{pn+1}
=\int_0^z\frac{dt}{1-t^p} \qquad \text{and} \qquad g_{p,\delta}(z)=(1+\delta)\sum_{n=1}^{\infty}
\frac{z^{pn+1}}{pn+1}
=(1+\delta)(h_p(z)-z).$$
Therefore, the Jacobian is
$$J_{f_{p,\delta}}(z)=|h_p'(z)|^2-|g_{p,\delta}'(z)|^2=\frac{1-(1+\delta)^2|z|^{2p}}
{|1-z^p|^2}.$$
If we choose $\,r_\delta=(1+\delta)^{-1/p},$ then $0<r_\delta<1$. Moreover at $z=r_\delta$, we have $J_{f_{p,\delta}}(r_\delta)=0.$ Hence  $J_{f_{p,\delta}}$ vanishes inside $\mathbb{D}$. Therefore $f$ is not univalent in $\mathbb{D}$. Since $\delta>0$ was arbitrary, the limiting constant $2$ cannot be increased. Therefore, the constant $2$ is sharp for $p=1,2.$

\medskip
\noindent
{\bf Remark 2.2}
\begin{itemize}
\item[(i)] Choosing $p=1$ in Theorem 2.1 gives the criterion for ordinary harmonic mappings \cite[Theorem 3.2]{Sahoo}.
\item[(ii)] Choosing $p=2$ in Theorem 2.1 gives the criterion for odd harmonic mappings \cite[Theorem 3.2]{Sahoo}.
\item[(iii)] Choosing $p=1$ and $g \equiv 0$ in Theorem 2.1 gives \cite[Theorem 1]{Mac}.
\item[(iv)] Choosing $p=2$ and $g \equiv 0$ in  Theorems 2.1 gives \cite[Theorem 2]{Mac}.
\end{itemize}

\medskip
The next result replaces monotonicity by a restriction on the total variation of weighted coefficient differences.

\medskip
\noindent
{\bf Theorem 2.3} Let $\{a_n\}_{n \ge 1}$ and $\{b_n\}_{n \ge 1}$ be two sequences of real numbers such that $a_1=1$, and $b_1=0$. For $n \ge 1$, suppose that $c_n=(n+1)a_{n+1}-na_n,\;\; d_n=(n+1)b_{n+1}-nb_n.$ If 
\begin{equation}\label{bound}
\qquad \qquad |c_1-1|+|d_1|+\sum_{n=1}^{\infty}
\left(|c_{n+1}-c_n|+
|d_{n+1}-d_n|\right)\le1,   
\end{equation}
then the harmonic mapping $\,f(z)=h(z)+\overline{g(z)}
=z+\sum_{n=2}^{\infty}a_nz^n+
\overline{\sum_{n=2}^{\infty}b_nz^n},
\; z\in\mathbb D,$ is sense-preserving, close-to-convex, and univalent in $\mathbb D$.

\medskip
\noindent
{\it Proof.} By hypothesis both sequences $\{c_n\}$ and $\{d_n\}$ are bounded. Since $\,na_n=1+\sum_{j=1}^{n-1}c_j,
\; nb_n=\sum_{j=1}^{n-1}d_j,$ the sequences $\{a_n\}$ and $\{b_n\}$ are bounded. Thus $h$ and $g$ are analytic in $\mathbb D$.
For $|\varepsilon|=1$, consider $\,F_\varepsilon(z)=h(z)+\varepsilon g(z).$ Now for $z \in \mathbb{D}$, we obtain
\begin{eqnarray*}
\text{Re}\{(1-z)^2F_\varepsilon'(z)\}
& \ge & 1-|c_1-1||z|-|d_1||z|-\sum_{n=1}^{\infty}
\bigl(|c_{n+1}-c_n|+
|d_{n+1}-d_n|\bigr)|z|^{n+1} \\
& \ge & 1-|z|>0
\end{eqnarray*}

The function $\,\psi(z)=z/(1-z)$ is convex in $\mathbb D$, and $\,\text{Re}
\left\{F_\varepsilon'(z)/\phi'(z)\right\}=\text{Re}\{(1-z)^2F_\varepsilon'(z)\}>0.$ Hence every $F_\varepsilon$ is close-to-convex in $\mathbb D$. Now using Lemma 1.1, $f=h+\overline g$ is close-to-convex and univalent in $\mathbb D$. 

Moreover, the last inequality implies $\,h'(z)+\varepsilon g'(z)\neq0$
for every $\varepsilon$ such that $|\varepsilon|=1$. Therefore, $|h'(z)|\neq|g'(z)|$ in $\mathbb D$.
Since $h'(0)=1$ and $g'(0)=0$, the continuity gives $\,|g'(z)|<|h'(z)|,\; z\in\mathbb D.$ Hence, $f$ is sense-preserving. This completes the proof.

\medskip
\noindent
{\bf Remark 2.4.} Theorem 2.3 remains applicable when the sequences $\{c_n\}$ and $\{d_n\}$ are nonmonotone. This distinguishes the bounded-variation method from coefficient criteria requiring monotonicity.

\medskip
The next Corollary relates Theorem 2.3 to three results of Sahoo and Wankhede \cite{Sahoo}.

\medskip
\noindent
{\bf Corollary 2.5.} Let $\{a_n\}_{n \ge 1}$ and $\{b_n\}_{n \ge 1}$ be two sequences of real numbers such that $a_1=1$, and $b_1=0$. Suppose that $c_n=(n+1)a_{n+1}-na_n,\;\; d_n=(n+1)b_{n+1}-nb_n,$ and  $\,L=\lim_{n\to\infty}(c_n+d_n).$ Then each of the following assumptions implies \eqref{bound}:
\begin{itemize}
\item[(i)] The sequences $\{c_n\}$ and $\{d_n\}$ are non-decreasing, $c_1\ge1$, $d_1\ge0$, and $L\le2$.
\item[(ii)] The sequences $\{c_n\}$ and $\{d_n\}$ are non-decreasing, $c_1\le1$, $d_1\ge0$, and $L\le2c_1$.
\item[(iii)] The sequences $\{c_n\}$ and $\{d_n\}$ are non-increasing, $0\le c_1\le1$, $d_1\ge0$, and $L\ge2d_1$.
\end{itemize}

\medskip
\noindent
{\it Proof.} Let $S$ denote the left-hand side of \eqref{bound}, and let $S_N=|c_1-1|+|d_1|+\sum_{n=1}^{N}
\left(|c_{n+1}-c_n|+
|d_{n+1}-d_n|\right)$. Under (i), we get
$$S=\lim_{N \to \infty} S_N=\lim_{N \to \infty} [c_N+d_N-1]=L-1\le1.$$
Under (ii), we get
$$S=\lim_{N \to \infty} S_N=\lim_{N \to \infty} [1-2c_1+c_N+d_N]=1-2c_1+L\le1.$$
Finally, under (iii), we get
$$S=\lim_{N \to \infty} S_N=\lim_{N \to \infty} [1+2d_1-c_N-d_N]=1+2d_1-L\le1.$$
Hence the result follows from Theorem 2.3.

\medskip
\noindent
{\bf Remark 2.6.}
Theorem 2.3 contains Theorems 3.4, 3.5, and 3.6 of Sahoo and Wankhede \cite{Sahoo}.

\section{Applications}

\medskip
\noindent
{\bf Theorem 3.1}
Let $p\in\{1,2\}$, $a,b>0$, and $c=a+b$. Suppose that $pab\ge1$ and $a+b\le(p+1)ab.$ Let $\lambda>0$ satisfy $\,p(\lambda+p(a+b))\Gamma(a+b)
\le2\lambda\Gamma(a)\Gamma(b).$ Then the harmonic mapping
\begin{equation}\label{appl1}
\mathcal H_{p,\lambda}^{a,b}(z)
={} z\,{}_2F_1(a,b;a+b;z^p)+\overline{
\frac{z^2}{\lambda}
{}_2F'_1(a,b;a+b+1;z^p)}    
\end{equation}
is close-to-convex and univalent in $\mathbb D$.

\medskip
\noindent
{\it Proof.}
We shall use Theorem 2.1. For $c=a+b$, and for $p\in\{1,2\}$, we write $\,\mathcal H_{p,\lambda}^{a,b}(z)=h_p(z)+\overline{g_p(z)}.$ Then
$$
h_p(z)=z+\sum_{n=1}^{\infty}A_{pn+1}z^{pn+1} \qquad \text{and} \qquad g_p(z)=\sum_{n=1}^{\infty}B_{pn+1}z^{pn+1}$$
where
\begin{equation}\label{appl2}
A_{pn+1}=\frac{(a)_n(b)_n}{(c)_n\,n!},
\qquad
B_{pn+1}=\frac{pn}{\lambda}
\frac{(a)_n(b)_n}{(c+1)_n\,n!}.    
\end{equation}
Also, $A_1=1$ and $B_1=0$.

Note that $A_{pn+1}\ge 0$ for all $n \ge 0$. Define $\,C_n=(pn+1)A_{pn+1},\; n\ge0.$ Then using $c=a+b$, we get for $n\ge0$,
\begin{eqnarray}\label{appl3}
C_{n+1}-C_n
&=&{}\frac{(p(n+1)+1)(a+n)(b+n)-(pn+1)(c+n)(n+1)}{(c+n)(n+1)}A_{pn+1}\\
&=& \frac{(pab-1)n+(p+1)ab-a-b}{(c+n)(n+1)}A_{pn+1} \ge 0.\notag
\end{eqnarray}
Hence $\,\{(pn+1)A_{pn+1}\}_{n\ge0}$ is non-decreasing.

Moreover, define $\,D_n=(pn+1)B_{pn+1},\; n\ge0,$ where $D_0=0$. From \eqref{appl2}, we have
$$B_{pn+1}=\frac{pcn}{\lambda(c+n)}A_{pn+1},
\; (n\ge1) \qquad \text{and}\quad D_1=(p+1)B_{p+1}
=\frac{p(p+1)ab}{\lambda(c+1)}>0.$$
Now using $c=a+b,$ we obtain
\begin{eqnarray*}
D_{n+1}-D_n
&=&\frac{(p(n+1)+1)(a+n)(b+n)-(pn+1)n(c+n+1)}{n(c+n+1)}B_{pn+1}\\
&=& \frac{(pab+pc-1)n+(p+1)ab}{n(c+n+1)}B_{pn+1}>0
\end{eqnarray*}
for every $n\ge1$. Hence, the sequence $\,\{(pn+1)B_{pn+1}\}_{n\ge0}$ is non-decreasing.

To complete the proof, it remains to show that $\,\lim_{n\to\infty}
(pn+1)(A_{pn+1}+B_{pn+1})\le2.$ By using the standard asymptotic formula for the Gamma function
\begin{equation*}
\qquad \qquad \qquad \lim_{n\to\infty}
\frac{\Gamma(n+\alpha)}
{\Gamma(n+\beta)}n^{\beta-\alpha}=1    
\end{equation*}
and using $c=a+b$, we obtain
\[\begin{aligned}
\lim_{n\to\infty}(pn+1)A_{pn+1}
&=\lim_{n\to\infty}\frac{\Gamma(c)}{\Gamma(a)\Gamma(b)}
\frac{pn+1}{n}
\left(\frac{\Gamma(n+a)}{\Gamma(n+1)} n^{1-a}\right)\left(
\frac{\Gamma(n+b)}{\Gamma(n+c)}n^{c-b}\right)\,n^{a+b-c} \\
&= \frac{p \,\Gamma(c)}{\Gamma(a)\Gamma(b)},
\end{aligned}\]
and
$$\lim_{n\to\infty}(pn+1)B_{pn+1}=\frac{pc}{\lambda}
\lim_{n\to\infty}
\frac{n}{c+n}(pn+1)A_{pn+1}=\frac{p^2c\Gamma(c)}
{\lambda\Gamma(a)\Gamma(b)}.$$
Thus
$$\lim_{n\to\infty}
(pn+1)(A_{pn+1}+B_{pn+1})=
\frac{p\Gamma(c)}{\Gamma(a)\Gamma(b)}
\left(1+\frac{pc}{\lambda}\right) \le 2$$
for $c=a+b.$ This concludes the proof.

\medskip
\noindent
{\bf Remark 3.2.}
\begin{itemize}
\item[(i)] Taking $p=1$ in Theorem 3.1, we obtain \cite[Theorem 4.1]{Sahoo}.
\item[(ii)] Taking $p=2$ in Theorem 3.1, we obtain \cite[Theorem 4.3]{Sahoo}.
\end{itemize}

\medskip
\noindent
{\bf Theorem 3.3}
Let $a,b,c>0$ satisfy $c>a+b+1.$ Suppose that $\alpha,\beta>0$ and
$$\alpha+\beta\leq
\frac{\Gamma(c-a)\Gamma(c-b)}
{4\Gamma(c)\Gamma(c-a-b-1)
\left[ab+2(c-a-b-1)\right]}.$$
Then the harmonic mapping
$$f(z)=\frac{z}{1-z}
+\alpha z^2{}_2F_1(a,b;c;z)
+\overline{\beta z^2{}_2F_1(a,b;c;z)}$$
is sense-preserving, close-to-convex, and univalent in $\mathbb D$.

\medskip
\noindent
{\it Proof.}
Let $f=h+\overline g$,\, then $\,h(z)=z+\sum_{n=2}^{\infty}a_nz^n\,$ and $\,g(z)=\sum_{n=2}^{\infty}b_nz^n,$ where for $n \ge 2$ and $C_{n}=\dfrac{(a)_{n}(b)_{n}}{(c)_{n} (n)!}$, we have $\,a_n=1+\alpha C_{n-2}$ and $\,b_n=\beta C_{n-2}$. Define $\,u_1=0\,$ and $\,u_n=n C_{n-2},\; n\ge2.$ Then $\,n a_n=n+\alpha u_n$ and $nb_n=\beta u_n\;\;(n\ge2).$ Thus the associated sequences are
$$c_n=(n+1)a_{n+1}-na_n=1+\alpha(u_{n+1}-u_n)$$
and
$$d_n=(n+1)b_{n+1}-nb_n=\beta(u_{n+1}-u_n).$$
Since $u_2=2$, we have $\,c_1=1+2\alpha,\; d_1=2\beta.$ Hence $\,|c_1-1|+|d_1|=2(\alpha+\beta).$ Moreover,
$$c_{n+1}-c_n=\alpha(u_{n+2}-2u_{n+1}+u_n) \quad \text{and} \quad d_{n+1}-d_n=\beta(u_{n+2}-2u_{n+1}+u_n).$$

Let $S$ denote the left-hand side of \eqref{bound}, then
$$S=2(\alpha+\beta)+(\alpha+\beta)
\sum_{n=1}^{\infty}
|u_{n+2}-2u_{n+1}+u_n|.$$
Since $u_n\geq0$, $|u_{n+2}-2u_{n+1}+u_n|
\leq u_{n+2}+2u_{n+1}+u_n.$ Now summing over $n\geq1$, and using $u_1=0,\ u_2=2$, we obtain
$$\begin{aligned} \sum_{n=1}^{\infty}
|u_{n+2}-2u_{n+1}+u_n| &\leq
\sum_{n=1}^{\infty}u_{n+2}
+2\sum_{n=1}^{\infty}u_{n+1}
+\sum_{n=1}^{\infty}u_n\\
&=
4\sum_{n=2}^{\infty}u_n-2.
\end{aligned}$$
Therefore $\,S\leq4(\alpha+\beta)
\sum_{n=2}^{\infty}u_n.$ Now,
$$\sum_{n=2}^{\infty}u_n=\sum_{k=0}^{\infty}
(k+2)\frac{(a)_k(b)_k}{(c)_k\,k!}=2F(1)+F'(1).$$
Since $c>a+b+1$, using Gauss's summation formula, eq (2), and the hypothesis, we obtain
$$S\leq4(\alpha+\beta)
\frac{\Gamma(c)\Gamma(c-a-b-1)}
{\Gamma(c-a)\Gamma(c-b)}\left[2(c-a-b-1)+ab\right]\leq1,$$
Therefore, Theorem 2.3 implies that $f$ is sense-preserving, close-to-convex, and univalent in $\mathbb D$.

\medskip
If $\{c_n\}$ and $\{d_n\}$ are monotone, their total variations can be expressed through their values and limits. However, if monotonicity is absent, the same criterion remains meaningful, provided the total variations are finite and sufficiently small. As a simple illustration, we give the following result.

\medskip
\noindent
{\bf Theorem 3.4}
Let $m$ be a nonnegative integer and let
$$A_n^{(m)}=\begin{cases}
\dfrac{(-m)_n(m+1)_n}{(n!)^2},&0\le n\le m,\\
0,&n>m.
\end{cases}$$
Suppose that $\;u_1=0,\; u_n=nA_{n-2}^{(m)} \;\; (n\geq2),$ and $\,v_m
=2+\sum_{n=1}^{m+2}
|u_{n+2}-2u_{n+1}+u_n|.$ If $\alpha,\beta>0\,$ satisfy $\,\alpha+\beta\leq 1/v_m,$
then the harmonic mapping
$$f_m(z)=\frac{z}{1-z}
+\alpha z^2P_m(1-2z)
+\overline{\beta z^2P_m(1-2z)}$$
is sense-preserving, close-to-convex, and univalent in $\mathbb D$.

\medskip
\noindent
{\it Proof.} Using the terminating hypergeometric series,
$$P_m(1-2z)=\sum_{k=0}^{m}
\frac{(-m)_k(m+1)_k}{(k!)^2}z^k,$$
we write $\,f_m=h_m+\overline{g_m},$ where
$$h_m(z)=\frac{z}{1-z}+\alpha z^2P_m(1-2z),
\qquad g_m(z)=\beta z^2P_m(1-2z).$$

For $n\geq2$, the Taylor coefficients are $\,a_n=1+\alpha A_{n-2}^{(m)},$ and $b_n=\beta A_{n-2}^{(m)}.$ Hence, the sequences appearing in Theorem 2.3 are
$$c_n=(n+1)a_{n+1}-na_n=1+\alpha(u_{n+1}-u_n),$$
and
$$d_n=(n+1)b_{n+1}-nb_n=\beta(u_{n+1}-u_n).$$

Since $u_1=0$ and $u_2=2$, we have $c_1=1+2\alpha,\; d_1=2\beta.$ Furthermore,
$$c_{n+1}-c_n
=\alpha(u_{n+2}-2u_{n+1}+u_n)$$
and
$$d_{n+1}-d_n
=\beta(u_{n+2}-2u_{n+1}+u_n).$$
Therefore, the left-hand side of (2) becomes
$$\begin{aligned}
S &=|c_1-1|+|d_1|+\sum_{n=1}^{\infty}
\left(|c_{n+1}-c_n|+|d_{n+1}-d_n|\right)\\
&=(\alpha+\beta)\left(
2+\sum_{n=1}^{m+2}
|u_{n+2}-2u_{n+1}+u_n|
\right)=(\alpha+\beta) v_m.
\end{aligned}$$
By hypothesis, $S\leq1$. This proves the assertion.

\medskip
\noindent
{\bf Remark 3.5.}
\begin{itemize}
\item[(i)] Taking $m=0$ in Theorem 3.4, we obtain $u_2=2,\; u_n=0\;(n\ge3),$ hence $v_0=8.$ Therefore, $\alpha+\beta\le 1/8\;$ is sufficient for the harmonic mapping
$$f_0(z)=\frac{z}{1-z}
+\alpha z^2+\overline{\beta z^2}$$
to be sense-preserving, close-to-convex, and univalent.
\item[(ii)] Taking $m=1$ in Theorem 3.4, we obtain: If \(\alpha,\beta>0\) satisfy $\alpha+\beta\le 1/32,$ then the harmonic mapping
\[f(z)=\frac{z}{1-z}
+\alpha z^2(1-2z)
+\overline{\beta z^2(1-2z)}\]
is sense-preserving, close-to-convex, and univalent in \(\mathbb D\).
\end{itemize}

\end{document}